\documentclass{gtart_a}
\pdfoutput=1

\title{Highly connected manifolds with positive Ricci curvature}

\author{Charles P Boyer}
\givenname{Charles P}
\surname{Boyer}
\address{Department of Mathematics and Statistics\\
University of New Mexico\\\newline
Albuquerque, NM 87131\\USA}
\email{cboyer@math.unm.edu}
\urladdr{}

\author{Krzysztof Galicki}
\givenname{Krzysztof}
\surname{Galicki}
\email{galicki@math.unm.edu}
\urladdr{}

\volumenumber{10}
\issuenumber{}
\publicationyear{2006}
\papernumber{48}
\lognumber{0651}
\startpage{2219}
\endpage{2235}

\doi{}
\MR{}
\Zbl{}

\arxivreference{math.DG/0508189}

\keyword{Sasakian manifold}
\keyword{positive Ricci curvature}
\keyword{links}
\keyword{diffeomorphism type}
\keyword{highly connected}
\keyword{}
\keyword{}
\subject{primary}{msc2000}{53C25}
\subject{secondary}{msc2000}{57R55}

\proposed{Gang Tian}
\seconded{Leonid Polterovich, Tobias Colding}
\received{5 September 2005}
\revised{22 October 2006}
\accepted{13 October 2006}
\published{29 November 2006}
\publishedonline{29 November 2006}
\corresponding{}
\editor{CPR}
\version{}

\let\xysavmatrix\xymatrix
\def\xymatrix{\disablesubscriptcorrection\xysavmatrix}
\AtBeginDocument{\let\bar\wbar\let\tilde\wtilde\let\hat\what}

\makeatletter
\def\cnewtheorem#1[#2]#3{\newtheorem{#1}{#3}[section]
\expandafter\let\csname c@#1\endcsname\c@theorem}

\numberwithin{equation}{section}

\renewcommand{\theenumi}{\textup{(\roman{enumi})}}

\newtheorem{theorem}{Theorem}[section]
\newtheorem{Thm}{Theorem}
\cnewtheorem{proposition}[theorem]{Proposition}
\cnewtheorem{corollary}[theorem]{Corollary}
\cnewtheorem{lemma}[theorem]{Lemma}

\theoremstyle{remark}
\cnewtheorem{example}[theorem]{Example}
\cnewtheorem{definition}[theorem]{Definition}
\cnewtheorem{conjecture}[theorem]{Conjecture}
\newtheorem{remark}{Remark}[section]

\newtheorem*{acknowledge}{Acknowledgements}

\makeatother  
\makeautorefname{Thm}{Theorem}

\makeatletter
\newcommand{\low}{\@ifnextchar^{}{^{\vphantom x}}}
\newcommand{\high}{\@ifnextchar_{}{_{\vphantom I}}}
\makeatother
\newcommand{\rssymb}[2]{\newcommand{#1}{{\mathrmsl{#2}}}}
\newcommand{\calsymb}[2]{\newcommand{#1}{{\mathcal{#2}}}}
\newcommand{\bbsymb}[2]{\newcommand{#1}{{\mathbb{#2}}}}
\newcommand{\lieoper}[2]{\newcommand{#1}{\mathop
  {\mathfrak{#2}\null}\nolimits}}
\newcommand{\oper}[3][n]{\newcommand{#2}{\mathop
  {\mathrm{#3}\null}\ifx n#1\nolimits\else\limits\fi}}
\newcommand{\rsoper}[3][n]{\newcommand{#2}{\mathop
  {\mathrmsl{#3}\null}\ifx n#1\nolimits\else\limits\fi}}
\bbsymb\bbf{F} \bbsymb\bbk{K} \bbsymb\bbi{I} \bbsymb\bbl{L} \bbsymb\bbo{O}
\bbsymb\bbj{J}
\bbsymb\bby{Y}
\bbsymb\bbp{P}
\bbsymb\bba{A}
\calsymb\cA{A} \calsymb\cB{B} \calsymb\cC{C} \calsymb\cD{D} \calsymb\cE{E}
\calsymb\cF{F} \calsymb\cG{G} \calsymb\cH{H} \calsymb\cI{I} \calsymb\cJ{J}
\calsymb\cK{K} \calsymb\cL{L} \calsymb\cM{M} \calsymb\cN{N} \calsymb\cO{O}
\calsymb\cP{P} \calsymb\cQ{Q} \calsymb\cR{R} \calsymb\cS{S} \calsymb\cT{T}
\calsymb\cU{U} \calsymb\cV{V} \calsymb\cW{W} \calsymb\cX{X} \calsymb\cY{Y}
\calsymb\cZ{Z}
\renewcommand{\geq}{\geqslant} 
\oper\End{End} \oper\Hom{Hom}                    
\oper\Sym{Sym} \oper\Skew{Skew}
\oper\Aut{Aut}                                   
\oper\GL{GL} \oper\SL{SL}\oper\Symp{Sp}
\oper\CO{CO} \oper\On{O} \oper\SO{SO} \oper\Pin{Pin} \oper\Spin{Spin}
\oper\CU{CU} \oper\Un{U} \oper\SU{SU}
\rsoper\Diff{Diff} \rsoper\SDiff{SDiff}
\lieoper\der{der}                                
\lieoper\gl{gl} \lieoper\sgl{sl}\lieoper\symp{sp}
\lieoper\cu{cu} \lieoper\un{u}  \lieoper\su{su}
\rsoper\Vect{Vect} \rsoper\Ham{Ham}
\newcommand{\norm}[2][]{|\mkern-2mu|#2|\mkern-2mu|
  _{\lower1pt\hbox{${}_{#1}$}}}
\newcommand{\Norm}[2][]{\bigl|\mkern-3mu\bigr|#2\bigr|\mkern-3mu\bigr|
  _{\lower1pt\hbox{${}_{#1}$}}}

\rsoper\dimn{dim}                           
\rsoper\grad{grad}                          
\rsoper\kernel{ker}\rsoper\image{im}        
\rsoper\alt{alt}   \rsoper\sym{sym}         
\rsoper\Ad{Ad}     \rsoper\ad{ad}           
\rsoper\CoAd{CoAd} \rsoper\coad{coad}       
\rsoper\trace{tr}  \rsoper\trfree{tf}       
\rsoper\detm{det}                           
\rsoper\Vol{Vol}                            
\rsoper\divg{div}                           
\rsoper\sign{sign}                          
\rssymb\iden{id}                            
\rssymb\vol{vol}                            
\oper\Imag{Im}\oper\Real{Re}                
\newcommand{\sd}{{\raise1pt\hbox{$\scriptscriptstyle +$}}}
\newcommand{\asd}{{\raise1pt\hbox{$\scriptscriptstyle -$}}}
\newcommand{\sdasd}{{\raise1pt\hbox{$\scriptscriptstyle\pm$}}}
\newcommand{\asdsd}{{\raise1pt\hbox{$\scriptscriptstyle\mp$}}}

\rsoper\scal{scal}
\def\kahl/{k\"ahler}
\def\Kahl/{K{\"a}hler}
\newcommand{\bfw}{{\mathbf w}}
\newcommand{\bfx}{{\mathbf x}}

\newcommand{\bfz}{{\mathbf z}}
\newcommand{\bfa}{{\mathbf a}}

\def\bba{{\mathbb A}}

\def\bbc{{\mathbb C}}

\def\bbf{{\mathbb F}}

\def\bbi{{\mathbb I}}
\def\bbj{{\mathbb J}}
\def\bbk{{\mathbb K}}
\def\bbl{{\mathbb L}}

\def\bbo{{\mathbb O}}
\def\bbp{{\mathbb P}}

\def\bby{{\mathbb Y}}
\def\bbz{{\mathbb Z}}

\def\grb{\beta}

\def\grs{\sigma}
\def\grt{\tau}

\def\grS{\Sigma}

\def\calf{{\mathcal F}}

\def\Se{Sasakian-Einstein }
\def\la#1{\hbox to #1pc{\leftarrowfill}}
\def\ra#1{\hbox to #1pc{\rightarrowfill}}

\def\lcm{{\rm lcm}}

\def\rank{\hbox{rank~}}

\def\rad{\hbox{rad}~}

\def\fract#1#2{\raise4pt\hbox{$ #1 \atop #2 $}}
\def\BOne{{\mathchoice {\rm 1\mskip-4mu l} {\rm 1\mskip-4mu l}
                          {\rm 1\mskip-4.5mu l} {\rm 1\mskip-5mu l}}}

\begin{document}

\begin{asciiabstract}
We prove the existence of Sasakian metrics with positive Ricci
curvature on certain highly connected odd dimensional manifolds. In
particular, we show that manifolds homeomorphic to the 2k-fold
connected sum of S^{2n-1} x S^{2n} admit Sasakian metrics with
positive Ricci curvature for all k. Furthermore, a formula for
computing the diffeomorphism types is given and tables are presented
for dimensions 7 and 11.
\end{asciiabstract}

\begin{htmlabstract}
We prove the existence of Sasakian metrics with positive Ricci
curvature on certain highly connected odd dimensional manifolds. In
particular, we show that manifolds homeomorphic to the 2k&ndash;fold
connected sum of S<sup>2n-1</sup>&times; S<sup>2n</sup> admit Sasakian metrics with
positive Ricci curvature for all k. Furthermore, a formula for
computing the diffeomorphism types is given and tables are presented
for dimensions 7 and 11.
\end{htmlabstract}

\begin{abstract}
We prove the existence of Sasakian metrics with positive Ricci
curvature on certain highly connected odd dimensional manifolds. In
particular, we show that manifolds homeomorphic to the 2k--fold
connected sum of $S^{2n-1}\times S^{2n}$ admit Sasakian metrics with
positive Ricci curvature for all k. Furthermore, a formula for
computing the diffeomorphism types is given and tables are presented
for dimensions 7 and 11.
\end{abstract}

\maketitle

\section*{Introduction}

An important problem in global Riemannian geometry is that of
describing the class of manifolds that admit metrics of positive Ricci
curvature. The only known obstructions for obtaining such metrics come
from either the classical Myers theorem or the obstructions to the
existence of positive scalar curvature which is fairly well understood
(see the recent review article by Rosenberg and Stolz \cite{RoSt01}
for discussion and references).  Given the lack of obstructions it
seems most natural to develop techniques for proving the existence of
positive Ricci curvature metrics.  Over the years several methods for
doing so have appeared. These include symmetry methods, bundle
constructions, surgery theory, algebro-geometric techniques. We refer
the reader to recent papers (Boyer, Galicki and Nakamaye
\cite{BGN03a,BGN03b}, Grove and Ziller \cite{GrZi02} and
Schwachh{\"o}fer and Tuschmann \cite{ScTu04}) for a discussion of the
history and pertinent references. In \cite{BGN03a,BGN03b}, with M
Nakamaye, we introduced a method for proving the existence of positive
Ricci curvature on odd dimensional manifolds which relies on a
transverse version of Yau's famous proof of the Calabi conjecture.
The odd dimensional manifolds to which this method has been applied
are hypersurfaces of isolated singularities coming from weighted
homogeneous polynomials.  All such manifolds are what are sometimes
called ``highly connected''.  So far, with collaborators, we have
applied our methods successfully mainly to rational homology spheres
\cite{BG02,BGN02a,BG03p}, homotopy spheres \cite{BGN03b,BGK03,BGKT03},
and connected sums of $S^2\times S^3$ \cite{BGN03c,BGN02b,BG03}; (see
also Koll\'ar \cite{Kol04}).  A recent review of our method can be
found in \cite{BG05a}. 

The purpose of this note is to prove the existence of Sasakian metrics
of positive Ricci curvature on certain odd dimensional highly
connected smooth manifolds. Manifolds of dimension $2n$ or $2n+1$ that
are $n-1$ connected are often referred to as {\it highly connected}
manifolds. They are relatively tractable and there is a classification
of such manifolds by C\,T\,C Wall \cite{Wal62,Wal67} and his students
(Barden \cite{Bar65} and Wilkens \cite{Wil72}) as well as Crowley
\cite{Cro01}. Our first result concerns dimension $4n-1$, where we
prove the following:

\begin{Thm}\label{thmA}
Let $n\geq 2$ be an integer, then
for each positive integer $k$ there exist Sasakian metrics with positive Ricci curvature in
$D_n(k)$ of the $|bP_{4n}|$ oriented diffeomorphism classes of the $(4n-1)$--manifolds $2k\#
(S^{2n-1}\times S^{2n})$ that bound a parallelizable manifold, where the number $D_n(k)$ is determined
by the explicit formula
\ref{numdiffeo} given in the Appendix. In particular, $D_n(1)=|bP_{4n}|$ for all $n\geq 2$, so
$2\# (S^{2n-1}\times
S^{2n})$ admits Sasakian metrics of positive Ricci curvature in every oriented diffeomorphism class.
\end{Thm}

Here $k\#(M_1\times M_2)$ denotes the $k$--fold connected sum of the manifold $M_1\times
M_2$.
In the case of the connected sums of products of standard spheres, metrics of positive Ricci
curvature have been constructed previously by Sha and Yang \cite{ShYa91}. However, the
existence of such metrics for the exotic differential structures appears to be new. In dimension
$4n+1$ we prove a somewhat weaker result:

\begin{Thm}\label{thmB}
For each pair of positive integers $(n,k)$ there exists an oriented  $(2n-1)$--connected 
$(4n+1)$--manifold $K$ with $H_{2n}(K,\bbz)$ free of rank $k$ which admits a Sasakian
metric of positive Ricci curvature. Furthermore, $K$ is diffeomorphic to one of the manifolds
$$\#k(S^{2n}\times S^{2n+1}), \quad \#(k-1)(S^{2n}\times S^{2n+1})\#T, \quad
\#k(S^{2n}\times S^{2n+1})\# \grS^{4n+1},$$ where
$T=T_1(S^{2n+1})$ is the unit tangent bundle of $S^{2n+1}$, and
$\grS^{4n+1}$ is the Kervaire sphere. For $k=1$ the manifolds
$$S^{2n}\times S^{2n+1}, \qquad (S^{2n}\times S^{2n+1})\# \grS^{4n+1}, \qquad T$$
all admit Sasakian metrics with positive Ricci curvature.
If $n=1,3$ then $\#k(S^{2n}\times S^{2n+1})$ admits a Sasakian metric with positive Ricci
curvature for all $k$.
\end{Thm}

The result here for $n=1$ was given previously in \cite{BGN03a}.
For highly connected rational homology spheres we have:

\begin{Thm}\label{thmC}
Every $(2n-2)$--connected oriented $(4n-1)$--manifold that is the boundary of
a parallelizable manifold whose homology group $H_{2n-1}(K,\bbz)$
is isomorphic to $\bbz_3$ admits Sasakian metrics with positive
Ricci curvature. There are precisely $2|bP_{4n}|$ such smooth oriented
manifolds.
\end{Thm}

There are two distinct oriented topological manifolds in this theorem and they are
distinguished by their linking form in $H_{2n-1}(K,\bbz)\approx \bbz_3$. However, they are equivalent as 
non-oriented manifolds. Moreover, each oriented manifold is comprised of $|bP_{4n}|$ distinct oriented 
diffeomorphism types.

\section{Highly connected manifolds}

The most obvious subclass of highly connected manifolds are the homotopy spheres which
were studied in detail in the seminal paper of Kervaire and Milnor \cite{KeMil63}. We briefly
summarize their results. Kervaire and Milnor defined an
Abelian group $\Theta_n$ which consists of equivalence classes of
homotopy spheres of dimension $n$ that are equivalent under
oriented h--cobordism. By Smale's famous h--cobordism theorem
\cite{Sma62b} (see also Milnor \cite{Mil65}) this implies
equivalence under oriented diffeomorphism. The group operation on
$\Theta_n$ is connected sum. Now $\Theta_n$ has an important
subgroup $bP_{n+1}$ which consists of equivalence classes of those oriented
homotopy spheres which are the boundary of a parallelizable
manifold. It is the subgroup $bP_{2n}$ that is important for us in the present work.  Kervaire
and Milnor proved:

\begin{enumerate}
\item $bP_{2m+1}=0$. \item $bP_{4m}$ ($m\geq2$) is cyclic
of order
$2^{2m-2}(2^{2m-1}-1)$ times the numerator of $\bigl({4B_m\over
m}\bigr)$, where $B_m$ is the $m$-th Bernoulli number.  Thus, for
example
 $|bP_8|=28, |bP_{12}|=992, |bP_{16}|=8128, |bP_{20}|=130,816$.
\item $bP_{4m+2}$ is either $0$ or $\bbz_2$.
\end{enumerate}

Determining which $bP_{4m+2}$ is $\{0\}$ and which is $\bbz_2$ has
proven to be difficult in general, and is still not completely
understood. If $m\neq 2^i-1$ for any $i\geq 3$, then Browder
\cite{Bro69} proved that $bP_{4m+2}=\bbz_2.$  However, $bP_{4m+2}$ is
the identity for $m=1,3,7,15$ (Mahowald and Tangora \cite{MaTa67} and
Barratt, Jones and Mahowald \cite{BaJoMa84}).  See Lance \cite{La00} for a
recent survey of results in this area and complete references. The
answer is still unknown in the remaining cases. Using surgery Kervaire
was the first to show that there is an exotic sphere in dimension $9$.
His construction works in all dimensions of the form $4m+1$, but as
just discussed they are not always exotic.

In analogy with the Kervaire--Milnor group $bP_{2n}$, Durfee
\cite{Dur77} defined the group $BP_{2n}$.  Actually, he first defined
this as a semigroup in \cite{Dur71}, and later in \cite{Dur77} with
the same notation denoted the corresponding Grothendieck group. Thus,
we have:

\begin{definition}\label{BPgroup}
For $n\geq 3$ let $SBP_{2n}$ denote the semigroup of oriented
diffeomorphism classes of closed oriented $(n-2)$--connected $(2n-1)$--manifolds that
bound parallelizable manifolds, and let $BP_{2n}$ denote its Grothendieck completion.
\end{definition}

As with $bP_{2n}$ mulitiplication in $SBP_{2n}$ is the connected sum operation, and the standard sphere is a 
two-sided identity. Thus, $SBP_{2n}$ is a monoid. Furthermore, the  Kervaire and Milnor
group $bP_{2n}$ is a subgroup of $BP_{2n}$ again by Smale's h--cobordism theorem.
Unless otherwise stated we shall heretofore assume that $n\geq 3$, and that ``manifold'' will mean oriented 
manifold. We are mainly interested
in those highly connected manifolds that can be realized as links of isolated hypersurface
singularities defined by weighted homogeneous polynomials, so we have:

\begin{definition}\label{whpset}
We denote by $WHP_{2n-1}$ the set of oriented diffeomorphism classes of smooth closed oriented $(2n-1)$ 
dimensional manifolds that can be realized as the link of an isolated hypersurface singularity of a weighted 
homogeneous polynomial in $\bbc^{n+1}$.
\end{definition}

It is easy to see that elements of $WHP_{2n-1}$ enjoy some nice properties.

\begin{theorem}\label{whptopthm}
Let $M\in WHP_{2n-1}$. Then
\begin{enumerate}
\item $M$ is highly connected, that is, it is $(n-2)$--connected.
\item $M$ is the boundary of a compact $(n-1)$--connected parallelizable manifold $V$ of
dimension $2n$ with $H_n(V,\bbz)$ free.
\item If $n$ is even the $(n-1)$st Betti number $b_{n-1}(M)$ is even.
\item $M$ is a spin manifold.
\end{enumerate}
\end{theorem}

\begin{proof}
Parts (i), (ii) and (iv) are well known. By \cite{BG01b}  $M$
admits a Sasakian structure, so its odd Betti numbers are even up to the middle dimension by
a well known result of Fujitani and Blair--Goldberg (cf \cite{BG05}) which proves (iii). 
\end{proof}

From (i) and (ii) of \fullref{whptopthm} one sees that there is a map 
$$\Phi\co WHP_{2n-1}\ra{2.1} BP_{2n}$$ 
which is the composition of the inclusion $WHP_{2n-1}\hookrightarrow SBP_{2n}$ with the natural semigroup 
homomorphism $SBP_{2n}\ra{1.6} BP_{2n}$.
The image $\Phi(WHP_{2n-1})$ is a
subset of $BP_{2n}$, and by (iii) of \fullref{whptopthm} it is a proper subset at least when $n$ is even. This 
can be contrasted with $bP_{2n}$
which by a result of Brieskorn \cite{Bri66} satisfies $bP_{2n}\cap
\Phi(WHP_{2n-1})=bP_{2n}$ if $n\geq 3$. Notice, however, that
$\Phi(WHP_{2n-1})$ is not generally a submonoid.

We now discuss invariants that distinguish elements of $BP_{2n}$. First, by Poincar\'e duality
the only non-vanishing homology groups occur in dimension $0,n-1,n$ and $2n-1$.
Moreover, $H_n(K,\bbz)$ is free and $\rank H_n(K,\bbz)= \rank H_{n-1}(K,\bbz)$. Thus, our
first invariant is the rank of $H_{n-1}(K,\bbz)$, so we define
\begin{equation}\label{BPk}
BP_{2n}(k)= \{K\in BP_{2n} ~|~ \rank H_{n-1}(K,\bbz)=k\}.
\end{equation}
This provides $BP_{2n}$ with a grading, namely
\begin{equation}\label{grBPk}
BP_{2n}=\bigoplus_k BP_{2n}(k)
\end{equation}
which is compatible with multiplication in $BP_{2n}$ in the sense that
\begin{equation}
\times\co  BP_{2n}(k_1)\times BP_{2n}(k_2)\ra{2.3} BP_{2n}(k_1+k_2).
\end{equation}
Note that $BP_{2n}(0)$ is the submonoid of highly connected rational homology spheres, and
that $BP_{2n}(k)$ is a $bP_{2n}$--module. 

The remaining known invariants (Wall \cite{Wal67}, Durfee \cite{Dur71,Dur77}) are a linking form on the torsion
subgroup of
$H_{n-1}(K)$ and a quadratic invariant on the $2n$--manifold whose boundary is $K$.
The precise nature of these invariants depends on whether $n$ is even or odd. For the case
$n$ even Durfee \cite{Dur77} shows that for $n\geq 3$ and $n\neq 4,8$ there is an exact
sequence
\begin{equation}\label{BPsequence}
0\ra{1.8} bP_{2n}\ra{1.8} BP_{2n}\fract{\Psi}{\ra{1.8}} \bbz\oplus
KQ(\bbz)\ra{1.8}0,
\end{equation}
where $KQ(\bbz)$ denotes the Grothendieck group of regular bilinear form modules over
$\bbz$. Let us describe the map $\Psi$. The projection onto the first factor is just the rank of
$H_{n-1}(K)$ while the projection onto the second factor is Wall's quadratic form \cite{Wal67}
which is essentially the classical linking form $b$ on the torsion subgroup of $H_{n-1}(K)$.
Any two manifolds $K_1,K_2\in BP_{2n}$ such that $\Psi(K_1)=\Psi(K_2)$ differ by a
homotopy sphere, ie, there is $\grS\in bP_{2n}$ such that $K_2\approx K_1\#\grS$ where
$\approx$ means diffeomorphic. It is well known (Brieskorn \cite{Bri66}) that for $n$ even the elements
$\grS\in bP_{2n}$ are determined by the signature of $V$. This completes the diffeomorphism
classification for $n\neq 4,8$ even. The cases $n=4,8$ are more complicated
(Wilkens \cite{Wil72}, Crowley \cite{Cro01}). Now, in
addition to the group $H_{n-1}(K)$ and the linking form $b$, there is an obstruction cocycle
$\hat{\grb}\in H^n(K,\pi_{n-1}(SO))\approx H^n(K,\bbz)$. The tangent bundle of $K$
restricted to the $(n-1)$--skeleton is trivial and $\hat{\grb}$ gives the obstruction to triviality on
the $n$--skeleton. If the torsion subgroup of $H_{n-1}(K)$ has odd order, then up to
decomposability these are all the invariants. However, if the torsion subgroup of $H_{n-1}(K)$
has even order, things are even more complicated, and the analysis in \cite{Wil72} was not
complete. It was recently completed in \cite{Cro01}. The important point for us is that if the
torsion subgroup of $H_{n-1}(K)$ vanishes, $K$ is determined completely up to
diffeomorphism by the rank of $H_{n-1}(K)$.
Summarizing we have:

\begin{theorem}\label{highlyconn1}
Let $M$ be a highly connected manifold in $BP_{4n}$ such that \break $H_{2n-1}(M,\bbz)
=\bbz^k$.
Then $M$ is diffeomorphic to $k\#(S^{2n-1}\times S^{2n})\# \grS^{4n-1}$ for some
$\grS^{4n-1}\in bP_{4n}$.
\end{theorem}

Notice that by a well-known result of Fujitani and Blair--Goldberg
(cf \cite{BG05}) $k\#(S^{2n-1}\times S^{2n})\# \grS^{4n-1}$ can
admit a Sasakian structure only if $k$ is even.

For the case $n$ odd the diffeomorphism classification was obtained by Wall
\cite{Wal67}, but for
our purposes, the presentation in \cite{Dur71} is more convenient. Let $K\in BP_{2n}$ with
$K=\partial V$, where $V$ can be taken as $(n-1)$--connected and parallelizable. In this case
the key invariant is a $\bbz_2$--quadratic form
$$\psi\co H_n(V,\bbz)/2H_n(V,\bbz)\ra{1.3} \bbz_2$$
defined as follows: Let $X$ be an embedded $n$--sphere in $V$ that
represents a non-trivial homology class in $H_n(V,\bbz)$, and let
$[X]$ denote its image in $H_n(V,\bbz)/2H_n(V,\bbz)$. Then
$\psi([X])$ is the characteristic class in the kernel \break
$\hbox{ker}(\pi_{n-1}(SO(n)\ra{1.3}\pi_{n-1}(SO))\approx \bbz_2$
of the normal bundle of $X$. Let $\rad \psi$ be the radical of
$\psi$, ie, the subspace of the $\bbz_2$--vector space
$H_n(V,\bbz)/2H_n(V,\bbz)$ where $\psi$ is singular. Then Durfee
\cite{Dur71} (see also \cite{DuKa75}) proves:

\begin{theorem}\label{highlyconn2}
Let $K_i\in BP_{2n}$ for $i=1,2$ with $n\geq 3$ odd be boundaries of parallelizable
$(n-1)$--connected $2n$ manifolds $V_i$ with $\bbz_2$ quadratic forms $\psi_i$. Suppose
that $H_{n-1}(K_1,\bbz)\approx H_{n-1}(K_2,\bbz)$, then
\begin{enumerate}
\item if $n=3$ or $7$, then $K_1$ and $K_2$ are diffeomorphic;
\item if the torsion subgroups of $H_{n-1}(K_i,\bbz)$ have odd
order and  \break $\psi_i|\rad \psi_i\equiv 0$ for $i=1,2$, then
$K_1\approx K_2\#(c(\psi_1)+c(\psi_2))\grS$, where $c$ is the Arf
invariant and $\grS$ is the Kervaire sphere, ie, the generator
of $bP_{2n};$ \item if the torsion subgroups of
$H_{n-1}(K_i,\bbz)$ have odd order and \break $\psi_i|\rad
\psi_i\not\equiv 0$ for $i=1,2$, then $K_1\approx K_2\approx
K_2\#\grS$.
\end{enumerate}
\end{theorem}

It is convenient to define $WHP_{2n-1}(k)$ to be the subset of $WHP_{2n-1}$ such that $H_{n-1}$ has rank $k$. 
Then (iii) of
\fullref{whptopthm} implies $WHP_{4n-1}(2k+1)=\emptyset$, whereas we shall see
that $WHP_{4n-1}(2k)\neq \emptyset$ as well as $WHP_{4n+1}(k)\neq \emptyset$ for all
$k$.
Recently, in the case $n=3$, Koll\'ar \cite{Kol04b,Kol05} has discovered strong restrictions on
the torsion subgroups of $H_{2}(K,\bbz)$ in order that $K$ admit a Sasakian structure which
implies that $\Phi(WHP_{6}(0))$ is a proper subset of $BP_{6}(0)$.  One certainly expects these
types of restrictions to persist in higher dimension as well.

\section{Branched covers and periodicity}

In this section we discuss some results of Durfee
and Kauffman \cite{DuKa75} concerning the periodicity of branched covers.
Let $K\subset S^{2n+1}$ be a simple fibered knot or link ($n\geq
1$), by which we mean an $(n-2)$ connected $(2n-1)$ embedded submanifold of $S^{2n+1}$
for which the Milnor fibration theorem holds. If $F$ is the Milnor fiber of the fibration $\phi\co 
S^{2n+1}-K\ra{1.5} S^1$ then the {\it monodromy map} $h\co H_n(F)\ra{1.3} H_n(F)$ is a
fundamental invariant of the link $K$. Let $K_k$ be a $k$--fold cyclic branched cover of
$S^{2n+1}$ branched along $K$. Then  Durfee and Kauffman \cite{DuKa75} show that there
is an exact sequence
\begin{equation}\label{DKexactseq}
H_n(F)\fract{\BOne+h+\cdots +h^{k-1}}{\ra{1.8}} H_n(F)\ra{1.5} H_n(K_k)\ra{1.5}
0.
\end{equation}
So homologically $K_k$ is determined by the cokernel of the map $\BOne +h+\cdots
+h^{k-1}$.
Now suppose that $K$ is a rational homology sphere and that the monodromy map $h$ of
$K$ has period $d$. Then since $\BOne -h$ is invertible, $\BOne +h+\cdots +h^{d-1}$ is the
zero
map in \ref{DKexactseq}, and this determines the homology of $K_d$. Summarizing
we have:

\begin{lemma}[Durfee--Kauffman]\label{DKlemma}
Let $K$ be a fibered knot in $S^{2n+1}$ which is a rational homology sphere such that the
monodromy map has period $d$.  Suppose further that $K_k$ is a $k$--fold cyclic cover of
$S^{2n+1}$ branched along $K$. Then
\begin{enumerate}
\item $H_n(K_d)\approx H_n(F)\approx \bbz^\mu$ where $\mu$ is the Milnor number of
$K$.
\item $H_*(K_{k+d})\approx H_*(K_k)$ for all $k>0$.
\item $H_*(K_{d-k})\approx H_*(K_k)$ for all $0<k<d$.
\end{enumerate}
\end{lemma}

\noindent Notice that (i) determines a large class of $n-1$ connected $2n+1$--manifolds
whose middle homology group $H_n$ is free, and in certain cases this determines the
manifold up to homeomorphism. Items (ii) and (iii) give a homological periodicity.

Durfee and Kauffman also show that there are both homeomorphism and diffeomorphism
periodicities in the case that $n$ is odd and $n\neq 1,3,7$. In particular in this case, when the
link $K$ is a rational homology sphere whose monodromy map has period $d$, $K_{k+d}$ is homeomorphic
to $K_k$. To obtain the
diffeomorphism periodicity let $\grs_k$ denote the signature of the intersection form on the
Milnor fiber $F_k$. Again assuming that $K$ is a rational homology sphere and $h$ has
periodicity $d$, one finds that $K_{k+d}$ is diffeomorphic to ${\grs_{d+1}\over 8}\grS\# K_k$
where ${\grs_{d+1}\over 8}\grS$ denotes ${\grs_{d+1}\over 8}$ copies of the Milnor sphere
$\grS$.  Here we state the slightly more general theorem of Durfee \cite[Theorem 6.4]{Dur77}:

\begin{theorem}\label{Dur6.4}
For even $n\neq 2,4,8$ let $K_i$ be $(n-2)$--connected manifolds that bound
parallelizable manifolds $V_i$, with $i=1,2$. Suppose that the
quadratic forms of $K_i$ are isomorphic and
$H_{n-1}(K_1,\bbz)\approx H_{n-1}(K_2,\bbz)$. Then
$\grs(V_2)-\grs(V_1)$ is divisible by $8$, and $K_2$ is
diffeomorphic to $K_1\#\frac{1}{8}(\grs(V_2)-\grs(V_1))\grS$ where
$\grs(V)$ is the Hirzebruch signature of $V$.
\end{theorem}

\begin{remark}\label{n=4,8rem}
\cite[Theorem 6.4]{Dur77} as well as \cite[Theorem 5.3]{DuKa75}
exclude the cases $n=4$ and $8$. However, it follows from \cite{Wil72} and \cite{Cro01} that
the diffeomorphism classification still holds in these cases since the links we are considering
here have no element
of even order in the torsion subgroup of $H_{n-1}$ (In fact the torsion subgroup vanishes in
the case above). This remark also pertains to the discussion for Theorem 3 below.
\end{remark}

\section{Positive Ricci curvature on links}\label{pos}

Recall \cite{BGN03a} that a Sasakian structure $(\xi,\eta,\Phi,g)$
is {\it positive} if the basic Chern class $c_1(\calf_\xi)$ of the
characteristic foliation $\calf_\xi$ is positive. The importance
of positive Sasakian structures comes from Theorem A of
\cite{BGN03a} which states that they give rise to Sasakian metrics
with positive Ricci curvature. An important ingredient in the
proof of this result is the `transverse Yau theorem' of El
Kacimi-Alaoui \cite{ElK}, or equivalently for the cases at hand, the orbifold version of Yau's theorem.
Now there is a natural induced Sasakian
structure on the link of a hypersurface singularity of a weighted
homogeneous polynomial \cite{BG01b}.  Combining this with
`orbifold adjunction theory' \cite{BG05} we obtain:

\begin{theorem}\label{posRiccilink}
Let $L_f$ be the link of an isolated hypersurface singularity of a weighted homogeneous
polynomial $f$ of degree $d$ and weight vector $\bfw$.  Suppose further that $|\bfw|-d>0$.
Then $L_f$ admits a Sasakian metric with
positive Ricci curvature.
\end{theorem}

It is a simple task to construct positive Sasakian structures on links by increasing the
dimension.

\begin{proposition}\label{stabprop}
Let $L_{f'}$ be the link of a weighted homogeneous polynomial\break $f'(z_2,\cdots,z_n)$ in $n-1$
variables with weight vector $\bfw'$ and degree $d'$. Assume that the origin in $\bbc^{n-1}$
is the only singularity so that $L_{f'}$ is smooth. Consider the weighted homogeneous
polynomial
$$f=z_0^2+z_1^2+f'$$
of degree $d=\lcm(2,d')$.
Then the link $L_f$ admits a Sasakian structure with positive Ricci curvature and
$b_{n-1}(L_f)=b_{n-3}(L_{f'})$.
\end{proposition}

\begin{proof}
There are two cases. If $d'$ is odd then the weight vector of $f$ is $\bfw=(d',d',2\bfw')$,
whereas, if $d'$ is even, then $\bfw=(\frac{d'}{2},\frac{d'}{2},\bfw')$. In the first case we have
$|\bfw|-d= d'+d'+2|\bfw'|-2d'=2|\bfw'|>0$, while in the second case $|\bfw|-d= \frac{d'}{2}
+\frac{d'}{2}+|\bfw'|-d'=|\bfw'|>0$. In either case $L_f$ admits a Sasakian metric with positive
Ricci curvature by \fullref{posRiccilink}. The equality of Betti numbers is well known and follows from 
a theorem of Sebastiani and Thom \cite{SeTh71,KaNe77}. 
\end{proof}

We note that it is easy to see that the appearance of the two $2$'s in $f$ implies that the klt
conditions used to imply the existence of \Se metrics \cite{BGK03,BG03p} cannot be satisfied.
So we can say nothing at present about the existence of \Se metrics on these links.

\section{Proofs of Theorems \ref{thmA}, \ref{thmB} and \ref{thmC}}\label{proofs}

The links that we need to prove Theorems \ref{thmA}--\ref{thmC} involve Brieskorn--Pham polynomials of the
form
\begin{equation}\label{BPpoly}
f_{p,q}=z_0^p+z_1^q+z_2^2+\cdots +z^2_n.
\end{equation}
The link associated with $f_{p,q}$ is
$$L_{p,q}=\{f_{p,q}=0\}\cap S^{2n+1}.$$
By \fullref{stabprop} all such links admit Sasakian metrics with positive Ricci curvature.
One can view $L_{p,q}$ as a $p$--fold branched cover of $S^{2n-1}$ branched over the link
$L_q$ defined by the polynomial
$$f_q=z_1^q+z_2^2+\cdots +z_n^2.$$

\begin{proof}[Proof of Theorem 1]
Here we need the link $L_{2(2k+1),2k+1}$, ie $p=2(2k+1),q=2k+1$, with $n$ even (here $n$ corresponds
to $2n$ in the statement of the theorem).
In this case the degree of $L_{2(2k+1),2k+1}$ is $d=2(2k+1)$ which is the period of the
monodromy map of the link $L_{2k+1}$. Furthermore,  $L_{2k+1}$ is a homotopy sphere
by the Brieskorn Graph Theorem \cite{Bri66} or \cite{BG05}.
Now the link $L_{2(2k+1),2k+1}$ is a $2(2k+1)$ branched cover of $S^{2n-1}$ branched
over $L_{2k+1}$, so by item (i) of \fullref{DKlemma}, we have
\begin{equation}\label{Lpqhom}
H_{n-1}(L_{2(2k+1),2k+1},\bbz)\approx H_n(L_{2(2k+1),2k+1},\bbz)\approx \bbz^\mu
=\bbz^{2k}.
\end{equation}
Here $\mu$ is the Milnor number \cite{Mil68} of the link $L_{2k+1}$ which is easily computed
by the
formula for Brieskorn polynomials, namely
$$\mu =\prod_{i=1}^n(a_i-1)=(2k+1-1)\cdot 1\cdots 1=2k.$$
\begin{remark}
Notice that the link $L_{2(2k+1),2k+1}$ can be obtained by iterating
\fullref{stabprop} beginning with the Brieskorn manifold
$M(2(2k+1),2k+1,2)$ which is described in \cite[Example 1 page
320]{Mil75}. As discussed there it is the total space of the circle
bundle with Chern number $-1$ over a Riemann surface of genus $k$.
\end{remark}

It now follows from \fullref{highlyconn1} that $L_{2(2k+1),2k+1}$ is diffeomorphic to
$2k\#(S^{n-1}\times S^n)\#\grS^{4n-1}$ for some $\grS^{4n-1}\in bP_{4n}$. (Here $n$ is as in the
statement of the theorem.)
We now use the periodicity results  of Durfee and Kauffman to determine the diffeomorphism
type. First we notice that \fullref{highlyconn1} together with 
\cite[Theorem 4.5]{DuKa75} imply that for
every positive integer $i$ and every positive integer $k$, the link $L_{2i(2k+1),2k+1}$ is
homeomorphic to the connected sum $2k\#(S^{n-1}\times S^n)$.
The diffeomorphism types are determined by \fullref{Dur6.4} (\cite[Theorem 6.4]{Dur77}, see also \cite[Theorem 5.3]{DuKa75}) together with \fullref{n=4,8rem}.
Let $F_{i,k}$ denote the Milnor fibre of the link $L_{2i(2k+1),2k+1}$ and $\grs(F_{i,k})$ its
Hirzebruch signature.  Then \fullref{Dur6.4} says that for each pair of positive integers
$i,j$ there is a diffeomorphism
\begin{equation}\label{diffeoit}
L_{2i(2k+1),2k+1}\approx\Bigl({\grs(F_{i,k})-\grs(F_{j,k})\over
8}\grS\Bigr)\# L_{2j(2k+1),2k+1},
\end{equation}
where $l\grS$ denotes the connected sum of $l$ copies of the Milnor sphere, and a minus sign
corresponds to reversing orientation. Actually this formula follows from a signature periodicity result of
Neumann as stated in \cite[Theorem 5.2]{DuKa75}. From Durfee's theorem the difference in
signatures is always divisible by 8, so this expression makes sense. Equation \ref{diffeoit} can
be iterated; so it is enough to consider the case $i=2$ and $j=1$. In order to determine how
many distinct diffeomorphism types occur in \ref{diffeoit}, we need to compute the signature
of the Milnor fibres. This is done in Appendix \ref{sigsection}. It is interesting to note that not all
diffeomorphism types can be attained. This ends the proof of \fullref{thmA}.
\end{proof}

\begin{proof}[Proof of Theorem 2]
Now we have $n$ odd (corresponding to $2n+1$ in the statement of the theorem) and there are several
cases. First we take $p=2(2k+1),q=2k+1$ as in
the proof of Theorem 1. Again this leads to the link $L_{2(2k+1),2k+1}$ with free homology
satisfying Equation \ref{Lpqhom} except now $n$ is odd. Next we
consider $q=2k$ in Equation \ref{BPpoly}. The link $L_{2k}$ of the Brieskorn--Pham
polynomial $f_{2k}=z_1^{2k}+z_2^2+\cdots +z_n^2$ is a rational homology sphere by the
Brieskorn Graph Theorem. Furthermore, its monodromy map has period $2k$. Then choosing
$p=2k$ in Equation \ref{BPpoly} the link $L_{2k,2k}$ is $2k$--fold branched cover over
$S^{2n+1}$ branched over the rational homology sphere $L_{2k}$, so by item (i) of \fullref{DKlemma}, we have
$$H_{n-1}(L_{2k,2k},\bbz)\approx H_n(L_{2k,2k},\bbz)\approx \bbz^\mu =\bbz^{2k-1}.$$
These two cases now give links whose middle homology groups are
free of arbitrary positive rank. However, unlike the case for $n$
even this does not determine the homeomorphism type unless $n=3,7$
in which case there is a unique diffeomorphism class. Indeed
\fullref{highlyconn2} implies we need to compute the quadratic
form $\psi$, and this appears to be quite difficult in all but the
simplest case. From \fullref{highlyconn2} one can conclude
\cite{Dur71} that if $M\in BP_{4n+2}$ with $H_{2n}(M,\bbz)$ free of
rank one, then it is homeomorphic to $S^{2n}\times S^{2n+1}$ or
the unit tangent bundle $T=T_1(S^{2n+1})$. (Now $n$ is as in the statement of the theorem). So the
diffeomorphism
types at most differ by an exotic Kervaire sphere $\grS^{4n+1}$.
Furthermore, $S^{2n}\times S^{2n+1},T$ and $(S^{2n}\times
S^{2n+1})\#\grS^{4n+1}$ generate the torsion-free submonoid of
$BP_{4n+2}$, there being relations in the monoid, namely,
$T\#T=2\#(S^{2n}\times S^{2n+1})$ and $T\# \grS^{4n+1}=T$ (Some
further relations may exist depending on $n$ such as
$T_1(S^3)\approx S^2\times S^3$). This proves the first statement
in \fullref{thmB}.

To prove the second statement we follow Durfee and Kauffman and consider a slightly
different Brieskorn--Pham polynomial, namely $z_0^{2k}+z_1^2+\cdots +z_n^2$. For $k=1$
we get as before a link $L_{2,2}$ whose middle homology group is free of rank one. Thus, it is
diffeomorphic to one of the three generators above by (i) of \fullref{DKlemma}. Now as
$k$ varies we have a homological periodicity by (ii) and (iii) of \fullref{DKlemma}. Durfee
and Kauffman show that there is an 8--fold diffeomorphism periodicity, and they compute the
$\psi$ invariant to show that
$$L_{2,2}\approx T, \quad L_{4,2}\approx (S^{2n}\times
S^{2n+1})\#\grS^{4n+1},$$ $$L_{6,2}\approx T\#\grS^{4n+1}\approx
T, \quad L_{8,2}\approx S^{2n}\times S^{2n+1}.$$ This proves
\fullref{thmB}.
\end{proof}

\begin{proof}[Proof of \fullref{thmC}]
This is essentially a corollary of  \cite[Proposition 7.2]{Dur77} where Durfee considers the link
$K_k$ of the Brieskorn--Pham polynomial $z_0^k+z_1^3+z_2^2+\cdots +z_n^2$ for even
$n\geq 4$. He shows that $H_n(K_2,\bbz)\approx H_n(K_4,\bbz)\approx \bbz_3$, but that
$K_2$ and $K_4$ have inequivalent linking forms. Furthermore, $K_{6l+2}$ is diffeomorphic
to $K_2\#(-1)^{\frac{n}{2}}l\grS^{4n-1}$ and $K_{6l+4}$ is diffeomorphic to
$K_4\#(-1)^{\frac{n}{2}}l\grS^{4n-1}$ where $\grS^{4n-1}=K_5$ is the Milnor generator.
\end{proof}

\appendix
\section{Computing the signature}\label{sigsection}

There are several known methods for computing the signature of the Milnor fibre $F$ of a
Brieskorn manifold in the case when $n$ is odd. This was first accomplished for homotopy
spheres by Brieskorn \cite{Bri66} and developed further by Hirzebruch and Zagier
\cite{Hir67,HiZa74}. Our
discussion follows that in \cite{Hir67}. Let $\bfa\in (\bbz^+)^{n+1}$ and write
$\bfa=(a_0,\cdots,a_n)$. Consider the Brieskorn manifold $M_\bfa$ defined by the link
$$\{z_0^{a_0}+\cdots +z_n^{a_n}=0\}\cap S^{2n+1}.$$
The Milnor fibre $F_\bfa$ can be represented by the Brieskorn manifold
$$\{\bfz\in \bbc^{n+1} ~|~z_0^{a_0}+\cdots +z_n^{a_n}=1\}.$$
For $n$ even the Hirzebruch signature of $F_\bfa$ is given by the function
$$\notag t(\bfa)= \#\{\bfx\in \bbz^{n+1} ~|~0<x_k<a_k~\hbox{and}~0<\sum_{j=0}^n{x_k\over
a_k} <1~\mod 2 \}
$$
\begin{equation}\label{Briesformula}
- \#\{\bfx\in \bbz^{n+1}
~|~0<x_k<a_k~\hbox{and}~1<\sum_{j=0}^n{x_k\over a_k} <2~\mod 2 \}.
\end{equation}
Using methods of Fourier analysis, Zagier has obtain the following formula for $t(\bfa):$
\begin{equation}\label{Zagformula}
t(\bfa)={(-1)^{{n\over 2}}\over N}\sum_{j=0}^{N-1}\cot{\pi(2j+1)\over
2N}\cot{\pi(2j+1)\over 2a_0}\cdots \cot{\pi(2j+1)\over 2a_n},
\end{equation}
where $N$ is any common multiple of the $a_i$'s.

We now adapt this formula to treat the link of the Brieskorn--Pham polynomial of Equation
\ref{BPpoly} with $N=2(2k+1)$, namely,
$\bfa=(2(2k+1),2k+1),2\cdots,2)$. Notice that we can always take the $N$ in Zagier's formula
\ref{Zagformula} to be the same as the $N$ in Equation \ref{BPpoly}. In this case we shall
denote $t(\bfa)$
by $t_d$ since the degree $d=2(2k+1)$ is the periodicity as well. Likewise, we denote by
$t_{2d}$ the signature $t(\bfa)$ with $\bfa=(4(2k+1),2k+1),2\cdots,2)$. We find
$$t_d={(-1)^{{n\over 2}}\over 4k+2}\sum_{j=0}^{4k+1}(-1)^j\cot^2{\pi(2j+1)\over
8k+4}\cot{\pi(2j+1)\over 4k+2},$$
and
$$t_{2d}={(-1)^{{n\over 2}}\over 8k+4}\sum_{j=0}^{8k+3}(-1)^j\cot^2{\pi(2j+1)\over
16k+8}\cot{\pi(2j+1)\over 4k+2}.$$ We want to compute
$\grt_k={|t_{2d}-t_d|\over 8}$. After some algebra we find that
$(64k+32)\grt_k$ equals
\begin{equation}\label{diffformula}
\sum_{j=0}^{8k+3}(-1)^j\cot{\pi(2j+1)\over 16k+8}
\Bigl(\cot{\pi(2j+1)\over 16k+8}-\cot{\pi(2j+1)\over
8k+4}\Bigr)\cot{\pi(2j+1)\over4k+2}.
\end{equation}
Now $\grt_k$ is always an integer, and by \ref{diffformula}  it is independent of
$n$.  We now define
\begin{equation}\label{numdiffeo}
D_n(k)={|bP_{4n}|\over \gcd(\grt_k,|bP_{4n}|)}.
\end{equation}
\def\strut{\vrule width 0pt depth 6pt height 11pt}
{\small
\begin{table}[t]
\centerline{
\vbox{\tabskip=0pt \offinterlineskip
\def\tablerule{\noalign{\hrule}}
\halign to310pt {\strut#& \vrule#\tabskip=1em plus2em&
     \hfil#& \vrule#& \hfil#& \vrule#& \hfil#& \vrule#&
     \hfil#& \vrule#\tabskip=0pt\cr\tablerule
\omit&height2pt&\multispan{7}&\cr &&\multispan{7}\hfil {\bf Table
1}\quad $2k\#(S^3\times S^4)$\hfil&\cr\tablerule &&\omit\hidewidth
$k$\hidewidth&& \omit\hidewidth $\grt_k$\hidewidth&&
\omit\hidewidth $D_2(k)$\hidewidth&& \omit\hidewidth ${D_2(k)\over
|bP_{8}|}$\hidewidth&\cr\tablerule
&&$1$&&$1$&&28&&$1$&\cr\tablerule
&&$2$&&$3$&&28&&$1$&\cr\tablerule &&$3$&&$6$&&14&&${1\over
2}$&\cr\tablerule &&$4$&&$10$&&14&&${1\over 2}$&\cr\tablerule
&&$5$&&$15$&&28&&1&\cr\tablerule &&$6$&&$21$&&4&&${1\over
7}$&\cr\tablerule &&$7$&&$28$&&1&&${1\over 28}$&\cr\tablerule
&&$8$&&$36$&&7&&${1\over 4}$&\cr\tablerule
&&$9$&&$45$&&28&&1&\cr\tablerule &&$10$&&$55$&&28&&1&\cr\tablerule
&&$20$&&$210$&&2&&${1\over 14}$&\cr\tablerule
&&$48$&&$1176$&&1&&${1\over 28}$&\cr\tablerule
&&$50$&&$1275$&&28&&1&\cr\tablerule &&$100$&&$5050$&&14&&${1\over
2}$&\cr\tablerule &&$496$&&$123256$&&1&&${1\over
28}$&\cr\tablerule &&$500$&&$125250$&&14&&${1\over
2}$&\cr\tablerule}} }
\end{table}
\begin{table}[t]
\centerline{ \vbox{\tabskip=0pt \offinterlineskip
\def\tablerule{\noalign{\hrule}}
\halign to310pt {\strut#& \vrule#\tabskip=1em plus2em&
     \hfil#& \vrule#& \hfil#& \vrule#& \hfil#& \vrule#&
     \hfil#& \vrule#\tabskip=0pt\cr\tablerule
\omit&height2pt&\multispan{7}&\cr &&\multispan{7}\hfil {\bf Table
2}\quad $2k\#(S^5\times S^6)$\hfil&\cr\tablerule &&\omit\hidewidth
$k$\hidewidth&& \omit\hidewidth $\grt_k$\hidewidth&&
\omit\hidewidth $D_3(k)$\hidewidth&& \omit\hidewidth ${D_3(k)\over
|bP_{12}|}$\hidewidth&\cr\tablerule
&&$1$&&$1$&&992&&1&\cr\tablerule
&&$2$&&$3$&&992&&1&\cr\tablerule
&&$3$&&$6$&&496&&${1\over 2}$&\cr\tablerule
&&$4$&&$10$&&496&&${1\over 2}$&\cr\tablerule
&&$5$&&$15$&&992&&1&\cr\tablerule
&&$6$&&$21$&&992&&1&\cr\tablerule
&&$7$&&$28$&&248&&${1\over 4}$&\cr\tablerule
&&$8$&&$36$&&248&&${1\over 4}$&\cr\tablerule
&&$9$&&$45$&&992&&1&\cr\tablerule
&&$10$&&$55$&&992&&1&\cr\tablerule
&&$31$&&$496$&&2&&${1\over 496}$&\cr\tablerule
&&$48$&&$1176$&&124&&${1\over 8}$&\cr\tablerule
&&$50$&&$1275$&&992&&1&\cr\tablerule
&&$62$&&$1953$&&32&&${1\over 31}$&\cr\tablerule
&&$124$&&$7750$&&16&&${1\over 62}$&\cr\tablerule
&&$248$&&$30876$&&8&&${1\over 124}$&\cr\tablerule
&&$496$&&$123256$&&4&&${1\over 248}$&\cr\tablerule
&&$500$&&$125250$&&496&&${1\over 2}$&\cr\tablerule
&&$992$&&$492528$&&2&&${1\over 496}$&\cr\tablerule}}
}\end{table}}
By Equation \ref{diffeoit}, $~D_n(k)$ represents the number of distinct
diffeomorphism types that can be represented by our construction. Using MAPLE we give two
tables consisting of a list of $\grt_k$ and $D_n(k)$ together with the ratio
$${D_2(k)\over |bP_{8}|}={1\over \gcd(\grt_k,|bP_{4n}|)}$$
for both the $7$--manifolds $\# 2k(S^3\times S^4)$ and the
$11$--manifolds $\# 2k(S^5\times S^6)$ for various values of $k$.

Notice that the prime factorization of $|bP_{4n}|$ consists of
high powers of two together with odd primes coming from the
Bernoulli numbers. Since $\grt_k$ is independent of $n$, this
gives rise to a bit of a pattern for the ratios ${D_n(k)\over
|bP_{4n}|}$. It is obvious that for $k=1$ all possible
diffeomorphism types occur, but this seems also to hold for $k=2$.
It is of course true whenever $|bP_{4n}|$ is relatively prime to
$3$. If we look at the next case namely, $bP_{16}$, we see that
$|bP_{16}|=8128=2^6\cdot 127$. Comparing this with
$|bP_{12}|=992=2^5\cdot 31$, we see that the same ratios will
occur for the case $\# 2k(S^7\times S^8)$ as for $\# 2k(S^5\times
S^6)$ for $k=1,\cdots,30$. It is interesting to contemplate whether the above gaps in the
diffeomorphism types occur as a consequence of our method or whether they indicate an
honest obstruction to the existence of positive Sasakian structures. At this stage we have no
way of knowing.

\begin{acknowledge}
We thank Stephan Stolz for helpful conversations, and David Wraith
for finding some typos and pointing out the need for some
clarifications. We also owe thanks to Walter Neumann and the Editorial Board of G\&T as well as the anonomous 
referee for further clarifications and corrections. The second author would like to thank
Max-Planck-Institut f\"ur Mathematik in Bonn for hospitality.
During the preparation of this work the authors
were partially supported by NSF grants DMS-0203219 and DMS-0504367.
\end{acknowledge}

\bibliographystyle{gtart}
\bibliography{link}

\end{document}